\documentclass{amsart}

\newtheorem{theorem}{Theorem}[section]
\newtheorem{lemma}[theorem]{Lemma}
\newtheorem{notations}[theorem]{Notations}
\newtheorem{proposition}[theorem]{Proposition}
\newtheorem{corollary}[theorem]{Corollary}

\theoremstyle{definition}
\newtheorem{definition}[theorem]{Definition}

\theoremstyle{remark}
\newtheorem{remark}[theorem]{Remark}

\numberwithin{equation}{section}

\title{On smooth surfaces in $\Pq$ containing a plane curve}
\author{Ph. Ellia - C. Folegatti}
\date{03/06/2004\\
Both authors supported by MIUR project "Geometry on Algebraic Varieties"}

\newcommand{\epf}{\ensuremath{\diamondsuit}}
\newcommand{\id}{\ensuremath{\mathcal{I}}}
\newcommand{\oc}{\ensuremath{\mathcal{O}}}
\newcommand{\fc}{\ensuremath{\mathcal{F}}}
\newcommand{\dc}{\ensuremath{\mathcal{D}}}

\newcommand{\nc}{\ensuremath{\mathcal{N}}}

\newcommand{\rc}{\ensuremath{\mathcal{R}}}

\newcommand{\pc}{\ensuremath{\mathcal{P}}}
\newcommand{\bc}{\ensuremath{\mathcal{B}}}
\newcommand{\Pt}{\mathbb{P}^3}
\newcommand{\Ptw}{\mathbb{P}^2}
\newcommand{\Pq}{\mathbb{P}^4}

\newcommand{\Pun}{\mathbb{P}^1}

\newcommand{\sig}{\Sigma}

\begin{document}
\maketitle

\section{Introduction}
We work over an algebraically closed field of characteristic zero.\\
In this paper we are dealing with smooth surfaces $S$ in $\Pq$ which contain a plane curve, $P$.\\
The first part contains some generalities about the linear system $|H-P|$, in particular we prove that its base locus has dimension zero and describe it. \\
In the second section we look at surfaces lying on a hypersurface of degree $s$ with a ($s$-2)-uple plane (we suppose $s \geq 4$), indeed if the surface does not lie on a hyperquadric, this implies that it contains a plane curve (lemma \ref{lem-plane}). The main results are the following.
\begin{theorem}
\label{th1}
Let $\sig \subset \Pq$ be an integral hypersurface of degree $s$ with a ($s$-2)-uple plane, then the degree of smooth surfaces $S \subset \sig$ with $q(S)=0$ is bounded by a function of $s$.
\end{theorem}
Then we restrict to the case of regular surfaces lying on a hyperquartic with singular locus of dimension two. It turns out that, if $deg(S) \geq 5$, the hyperquartic must have a double plane (lemma \ref{lem-planeQ}). In this situation we can compute an effective bound.
\begin{theorem}
\label{th2}
Let $S \subset \Pq$ be a smooth surface with $q(S)=0$ and lying on a quartic hypersurface $\sig$, such that $Sing(\sig)$ has dimension two, then $d=deg(S) \leq 40$.
\end{theorem}
The assumption $q(S)=0$ is due to technical reasons, in fact we believe that it is not strictly necessary (see \ref{general}).\\
Theorem \ref{th2} is of some interest for the classification of surfaces not of general type, since in this case one has to look only at surfaces lying on low degree hypersurfaces. For similar results concerning smooth surfaces on hyperquartics with isolated singularities see \cite{EF}.
\section{Smooth surfaces containing a plane curve}
Let $S \subset \Pq$ be a smooth, non degenerate surface, of degree $d$, containing a plane curve, $P$, of degree $p$. If $p \geq 2$, there is a unique plane, $\Pi$, containing $P$; otherwise if $P$ is a line, there are $\infty^2$ such planes, we just choose one of them and call it $\Pi$. We assume that $P$ is the one-dimensional part of $\Pi \cap S$. We denote by $\delta$ the linear system cut out on $S$, residually to $P$, by the hyperplanes containing $\Pi$. Since, by Severi's theorem, $H^0(\oc _S(1)) \simeq H^0(\oc _{\Pq}(1))$ (we assume $S$ is not a Veronese surface), $\delta = |H-P|$ if $p \geq 2$; if $P$ is a line, $\delta$ is a pencil in the $\infty ^2$ linear system $|H-P|$. Finally we will denote by $Y_H$ the element of $\delta$ cut out by the hyperplane $H$, and by $C_H = P \cup Y_H$, the corresponding hyperplane section of $S$.
\begin{lemma}
\label{lem1}
(i) The curve $P$ is reduced and the base locus of $\delta$ is empty or zero-dimensional and contained in $\Pi$. The general element $Y_H \in \delta$ is smooth out of $\Pi$ and doesn't have any component in $\Pi$.\\
(ii) If $p=1$, the linear system $|H-P|$ is base point free.
\end{lemma}
\textit{Proof:} (i) Clearly the base locus of $\delta$ is contained in $\Pi$. Assume an irreducible component of $P$, $P_1$, is in the base locus of $\delta$. Then, for every $H$ through $\Pi$, $C_H =H\cap S$ is singular along $P_1$. It follows that $T_xS \subset H$, for every $x \in P_1$. Since this holds for every $H$ through $\Pi$, we get $T_xS = \Pi$, $\forall x \in P_1$, but this contradicts Zak's theorem (\cite{Z}) which states that the Gauss map is finite. The same argument shows that $P$ is reduced. We conclude by Bertini's theorem.\\
(ii) Assume $P$ is a line. Clearly the base locus of $|H-P|$ is contained in $P$. Take $x \in P$. Now let $H$ be an hyperplane containing $P$ but not containing $T_xS$, then $C_H=P\cup Y_H$ is smooth at $x$, so $x \notin Y_H$.\epf
\begin{remark} (i) If $p=1$, $|H-P|$ is base point free and yields a morphism $f: S \to \Ptw$, which is nothing else than the projection from the line $P$. If there is no plane curve on $S$ in a plane through $P$, $f$ is a finite morphism of degree $d-2+P^2$.\\
(ii) Let $S \subset \Pq$ be an elliptic scroll, then $S$ contains a one dimensional family of cubic plane curves which are unisecants. If $P$ is such a cubic, and if $H$ is a general hyperplane through $P$, then $H \cap S=P \cup f \cup f'$, where $f,f'$ are two rulings. This shows that the general curve $Y_H \in |H-P|$ need not be irreducible.
\end{remark}
Since $\delta$ is a pencil and since the base locus, $\bc$, is zero-dimensional, the degree of $\bc$ is $(H-P)^2$. Now we give a geometric description of $\bc$. Let $Z:=\Pi \cap S$, $Z$ is a 1-dimensional subscheme of $\Pi$ (and also of $S$) and is composed by $P$ and possibly by some 0-dimensional component, which may be isolated or embedded in $P$.
\begin{definition}
We define $\rc$ as the residual scheme of $Z$ with respect to $P$, hence $\id _{\rc} = (\id _Z:\id _P)$.
\end{definition}
Since $\rc \subset Z$, we can view $\rc$ as a subscheme of $\Pi$ or of $S$.
\begin{lemma}
\label{lem-RB}
We have $\bc = \rc$.
\end{lemma}
\textit{Proof:} We observe that $\rc \subset \bc$ and that $deg(\bc)=d-2p+P^2$, then we only have to compute $deg(\rc)$.\\
Considering a section of $\omega _S(2)$ (which is always globally generated), we can associate to $S$ a reflexive sheaf $\fc$ of rank two and an exact sequence: $0 \to \oc_{\Pq} \stackrel{t}{\to} \fc \to \id_S(3) \to 0$ such that $(t)_0=S$. The singular locus of $\fc$ is a divisor in $|2H+K|$ and the Chern classes of $\fc$ are $c_1=3$, $c_2=d$.\\
We can restrict the sequence above to $\Pi$ and get a section $0 \to \oc_{\Pi} \stackrel{t_{\Pi}}{\to} \fc_{\Pi}$. Clearly $P \subset (t_{\Pi})_0$ then dividing by an equation of $P$ we get a non-zero section $\bar{t}_{\Pi}$ of $\fc_{\Pi}(-p))$. We compute $deg((\bar{t}_{\Pi})_0)= c_2(\fc_{\Pi}(-p))=c_2(\fc(-p))=-3p+d+p^2$.
The section $\bar{t}_{\Pi}$ will vanish on $\rc$ and on the intersection with $\Pi$ of the singular locus of $\fc$, which is a curve $X \in |2H+K|$. Thus $(\bar{t}_{\Pi})_0 = \rc \cup (X \cap \Pi)$. When we restrict to $\Pi$ we have $X \cap \Pi= X \cap P$ and we get $\sharp(X \cap \Pi)= (2H+K)P=2p+PK$.\\
It follows that $deg(\rc)= -5p+d+p^2-PK$. Now we use adjunction to get $PK=p^2-3p-P^2$ and combining with the previous equation we obtain the result.\epf
\begin{remark}
(a) There is a cheaper proof of this result. By looking through the lines of \cite{F}, page 155, we infer that $deg(\rc)=d-2p+P^2$.\\
Indeed we can see $S \cap \Pi$ as the intersection of two hyperplane divisors on $S$, $H_1$ and $H_2$ such that $H_1 \cap H_2= \Pi$. Moreover $P$ is a Weil divisor on the smooth surface $S$, hence a Cartier divisor. Then we compute the equivalence of $P$ in the intersection $H_1 \cap H_2$, namely $(H_1 \cdot H_2)^P=(H_1 +H_2 -P) \cdot P=2p-P^2$. This means that the "exceeding" curve $P$ counts for $2p-P^2$ points in $H_1 \cap H_2$, thus the degree of its zero-dimensional component, $\rc$, drops by $2p-P^2$. It follows that $\deg(\rc)=d-2p+P^2$, hence the result.\\
\end{remark}
\section{Degree $s$ hypersurfaces with a ($s$-2)-uple plane}
\begin{lemma}
\label{lem-plane}
If $S \subset \Pq$ is a smooth surface, lying on a degree $s$ integral hypersurface $\sig$ with a ($s$-2)-uple plane, then $S$ contains a plane curve or $h^0(\id_S(2)) \neq 0$.
\end{lemma}
\textit{Proof:} Let $\Pi$ be the($s$-2)-uple plane in $\sig$ and let $H$ be an hyperplane containing $\Pi$, then $H \cap \sig = (s-2)\Pi \cup Q$, where $Q$ is a quadric surface and $C_H=S \cap H \subset (s-2)\Pi \cup Q$. If $dim(C_H \cap \Pi)=0$, then $C_H \subset Q$, i.e. $h^0(\id_{C_H}(2)) \neq 0 $ and the same holds for $S$. Then we can assume $dim(C_H \cap \Pi)=1$ and this is equivalent to say that $S$ contains a plane curve.\epf
\begin{notations}
Let $\sig \subset \Pq$ be an integral hypersurface of degree $s$ containing a plane, $\Pi$, in its singular locus, with multiplicity $(s-2)$. Let $S \subset \sig$ be a smooth surface. If $h^0(\id _S(2)) \neq 0$, then $d:=deg(S) \leq 2s$. From now on we assume $h^0(\id _S(2))=0$. By Lemma \ref{lem-plane}, $dim(S \cap \Pi )=1$ and we denote by $P$ the 1-dimensional component of $\Pi \cap S$, also we let $p:=deg(P)$.\\
We assume $q(S)=0$, this assumption implies that every hyperplane section $C=H \cap S$ is linearly normal in $H \simeq \Pt$.\\
If $H$ is an hyperplane through $\Pi$, we denote by $C =Y_H \cup P$ the hyperplane section $H \cap S$. We have $C \subset \sig \cap H = (s-2)\Pi \cup Q_H$, where $Q_H$ is a quadric surface. By Lemma \ref{lem1}, if $H$ is general, $Y_H \subset Q_H$. If we restrict to $\Pi$, the $q_H=Q_H \cap \Pi$ form, as $H$ varies, a family of conics in $\Pi$. Let us set $\bc_q=\displaystyle{ \bigcap_{H \supset \Pi} q_H}$, $\bc_q$ is the base locus of the conics $q_H$. Since $Y_H \cap \Pi \subset Q_H \cap \Pi =q_H$, we have $\rc \subset \bc_q$.\\
Recall that if $\mu =c_2(\nc_S(-s))= d(d+s(s-4))-s(2\pi-2)$ ($\pi$ is the sectional genus of $S$), then by lemma 1 of \cite{EP}: $0 \leq \mu \leq (s-1)^2d-D(3H+K)$ where $D$ is the one dimensional part of the intersection of $S$ with $Sing(\sig)$. In our situation $P \subset D$, so $\mu \leq (s-1)^2d-P(3H+K)=(s-1)^2d-3p-PK$. By adjunction we compute $P^2+PK=p^2-3p$ and then $\mu \leq (s-1)^2d-p^2+P^2= s(s-2)d-p^2+2p+r$ (since $r=d-2p+P^2$).
\end{notations}
\begin{lemma}
\label{lem-multiple}
With the notations above, the base locus $\bc_q$ of the conics $q_H$ is ($s$-1)-uple for $\sig$.
\end{lemma}
\textit{Proof:} We assume the plane $\Pi$ is given by $x_0=x_1=0$, thus if $\phi=0$ is an equation of $\sig$ we have $\phi \in (x_0,x_1)^{s-2}$. We can write for example $\phi=\displaystyle{\sum_{i=0}^{s-2} Q_i(x_0,x_1,x_2,x_3,x_4) x_0^i x_1^{s-2-i}}$ where the $Q_i$ are quadratic forms.\\
The general hyperplane $H_{\alpha}$ containing $\Pi$ has an equation of the form $x_0=\alpha x_1$, $\alpha \in k$, we consider $\phi _{|H_{\alpha}}$, namely the equation of the surface $\sig \cap H_{\alpha}$:\\
$\phi_{|H_{\alpha}}=\displaystyle{\sum_{i=0}^{s-2} Q_i(\alpha x_1,x_1,x_2,x_3,x_4) \alpha^i x_1^{s-2}}=x_1^{s-2}\displaystyle{\sum_{i=0}^{s-2} Q_i(\alpha x_1,x_1,x_2,x_3,x_4) \alpha^i}$.\\
Clearly $\displaystyle{\sum_{i=0}^{s-2} Q_i(\alpha x_1,x_1,x_2,x_3,x_4) \alpha^i}=0$ is an equation defining $Q_H$ for the hyperplane $H_{\alpha}$. Let $x=(0:0:x_2:x_3:x_4)$ be a point in $\bc_q$, hence $\displaystyle{\sum_{i=0}^{s-2} Q_i(x) \alpha^i}=0$ for all $\alpha \in k$ and this implies that $Q_i(x)=0$.\\
Now if we look at the ($s$-2)-th derivatives of $\phi$, we see that they all vanish in a point $x \in \bc_q$, equivalently $x$ is a ($s$-1)-uple point for $\sig$.\epf
\begin{lemma}
\label{lem-planeQ}
If $S \subset \Pq$ is a smooth surface with $q(S)=0$, lying on a quartic hypersurface $\sig$ having singular locus of dimension two, then, if $deg(S) \geq 5$, the component of dimension two in $Sing(\sig)$ is a plane (or a union of planes) and $S$ contains a plane curve.
\end{lemma}
\textit{Proof:} Let us suppose that $Sing(\sig)$ contains an irreducible surface of degree $>1$, then the general hyperplane section $S \cap H = C$ lies on $F=\sig \cap H$, which is a quartic surface of $\Pt$ having an irreducible curve of degree $>1$ in its singular locus. From the classification of quartic surfaces in $\Pt$ it follows that such a surface is a projection of a quartic surface $F' \subset \Pq$, then $F$ is not linearly normal. Since $C$ is linearly normal and smooth, the curve $C' \subset F'$ projecting down to $C$ must be degenerate and thus $d=deg(C') \leq 4$. So we may assume that the singular locus of $\sig$ does not contain irreducible surfaces of degree $>1$.
Thus $Sing(\sig)$ contains a plane, say $\Pi$, which is double in $\sig$. Indeed $\sig$ cannot have a triple plane, otherwise $F=\sig \cap H$ would be a quartic surface in $\Pt$ with a triple line, and we argue as before because such a surface is not linearly normal in $\Pt$. By lemma \ref{lem-plane}, $S$ contains a plane curve.\epf\\
\\
\textit{Proof of theorems \ref{th1} and \ref{th2}}\\
We must distinguish between different cases, according to the behaviour of the curves $q_H$. Note that it is not possible that $q_H=0$ for every $H$; indeed if it were so, $\Pi$ would be ($s$-1)-uple for $\sig$. Then for all hyperplanes $H \supset \Pi$, $\sig \cap H=(s-1)\Pi \cup \Pi_H$, where $\Pi_H$ is a plane. With notations as above we could say that $Q_H=\Pi \cup \Pi_H$, but we know by lemma \ref{lem1} that, if $H$ is general, $Y_H$ does not have any component in $\Pi$, then $Y_H \subset \Pi_H$ is a plane curve and $h^0(\id_C(2)) \neq 0$: absurd.\\
So we are left with the following possibilities. The conics may move, i.e. vary as $H$ varies, so that at least two of them intersect properly, then $dim(\bc_q)=0$; conversely they may all be equal to a fixed conic $q$ or they can be all reducible and contain a fixed line $D$, while the remaining line is moving. Observe that there are always two possibilities: the one-dimensional part of $\bc_q$ could be contained in $S$ or not.
The starting point of the proof is trying to show that $h^1(\id_C(2))=0$ where $C=Y_H\cup P$. Indeed if it is so, then by $0 \to \id_S(1) \to \id_S(2) \to \id_C(2) \to 0$ we obtain $h^1(\id_S(2))=0$. Then using $0 \to \id_S(2) \to \id_S(3) \to \id_C(3) \to 0$ and the fact that $h^0(\id_C(3)) \neq 0$ we  get that $h^0(\id_S(3)) \neq 0$ and this implies $d \leq 3s$.\\
The proof will follow from the lemmas below.
\begin{lemma}
\label{small-genus}
If $p_a(Y_H) \leq 2(d-p-4)$ and if $r \leq 4$, then $d$ is bounded by a function of $s$. More precisely if $s=4$, $d \leq 40$.
\end{lemma}
\textit{Proof:} We have $\pi=p_a(Y_H)+ \frac{(p-1)(p-2)}{2}+d-p-r-1$, so $\pi-1 \leq 3(d-p)+\frac{p^2-3p}{2}-9-r$. Since $\mu \leq s(s-2)d-p^2+2p+r$ and on the other hand $\mu= d(d+s^2-4s)-2s(\pi-1)$, this yields: $\pi-1 \geq \frac{d^2-2sd+p^2-2p-r}{2s}$.\\
Now comparing the lower and the upper bound on $\pi-1$ we obtain: $ d^2-8sd+p^2(1-s)+p(9s-2)+18s+r(2s-1) \leq 0$ and since $r \geq 0$ it becomes: $ d^2-8sd+p^2(1-s)+p(9s-2)+18s \leq 0$. This implies $d \leq 4s + \sqrt{\Delta} \:\: (*)$, where $\Delta=16s^2 + p^2(s-1)-p(9s-2)-18s$. A short calculation shows that $\sqrt{\Delta} \leq p\sqrt{s-1} + 4s$ for all $s \geq 0$. In conclusion: $d \leq 8s + p \sqrt{s-1}$.\\
We take into account again the relation: $0 \leq \mu \leq s(s-2)d-p^2+2p+4$ and using the bound on $d$ stated above it becomes: $s(s-2)(8s+p\sqrt{s-1}) \geq p^2-2p-4$. This implies that $p$ is bounded by a function of $s$. We conclude since $d \leq 8s+p \sqrt{s-1}$.
\par

If $s=4$ we give a better bound for $\sqrt{\Delta}$, indeed $\sqrt{\Delta} \leq p\sqrt{3}-8$ if $p \geq 19$, thus $d \leq 8 + p\sqrt{3}$. The same relation used above now gives: $8d \geq p^2-2p-4$, hence $p^2-2p-8\sqrt{3}p-68 \leq 0$, which implies $p \leq 19$ and consequently by $(*)$: $d \leq 40$. On the other hand if $p \leq 18$, again by $(*)$ we have $d \leq 39$.\epf
\begin{lemma}
\label{rleq4}
If $r \leq 4$ and if $\rc$ does not contain three collinear points, then $d$ is bounded by a function of $s$. In particular if $s=4$, $d \leq 40$.
\end{lemma}
\textit{Proof:} Assume first $Q_H$ is a smooth quadric surface. We have $Y_H \cap \Pi=Y_H \cap P + \rc$, so $0 \to \id_C(2) \to \id_P(2) \to \oc_{Y_H}(\rc+1) \to 0$. The curve $Y_H$ has bidegree $(a,b)$, $a \leq b$. We may assume $a \geq 4$, otherwise $p_a(Y_H) \leq 2(d-p-4)$ and we conclude by lemma \ref{small-genus}.\\
Thus $Y_H$ is linearly normal. We have $h^0(\oc_{Y_H}(1+\rc))=4$ if and only if $\rc$ gives independent conditions to $\omega_{Y_H}(-1)$. This is equivalent to say that $\rc$ gives independent conditions to the curves of bidegree $(a-3,b-3)$. If $a=b=4$, then $deg(Y_H)=d-p=8$ and using $s(s-2)d-p^2+2p+4 \geq 0$ we get $0 \leq -d^2+d(18+s(s-2))-76$. This shows that $d$ is bounded by a function of $s$, in particular if $s=4$, $d \leq 22$. So we may assume $a \geq 4$, $b \geq 5$ and since $r \leq 4$ and no three points of $\rc$ are collinear, the curves of bidegree $(a-3,b-3)$ separate the points of $\rc$. It follows that the map $H^0(\id_P(2)) \to H^0(\oc_{Y_H}(1+\rc))$ is surjective, hence $h^1(\id_C(2))=0$. As said before, this implies $d \leq 3s$.\\
Now we suppose $Q_H$ is an irreducible quadric cone (recall that every reduced curve on a quadric cone is a.C.M.). If $d-p$ is even, then $Y_H$ is a complete intersection $(\frac{d-p}{2},2)$ and $\omega_{Y_H} \cong \oc_{Y_H}(\frac{d-p}{2}-2)$. So if $\frac{d-p}{2}-3 \geq 3$, arguing as above, we get $h^0(\oc_{Y_H}(1+\rc))=4$. On the other hand if this condition is not satisfied then $d-p \leq 11$, i.e. $p \geq d-11$. Recall that $0 \leq \mu \leq s(s-2)d-p^2+2p+4$; it follows that $(d-11)(d-13) \leq s(s-2)d+4$ and, for fixed $s$, this implies that $d$ is bounded. If $s=4$ we have: $d^2-32d+139 \leq 0$ which yields $d \leq 26$.\\
If $d-p$ is odd, $Y_H$ is linked to a line $L$ by a complete intersection $T$ of type $(\frac{d-p+1}{2},2)$. Since $L$ can be any ruling of $Q_H$, we may assume $L \cap \rc = \emptyset$. The exact sequence of liaison: $0 \to \id_{T}(\frac{d-p-5}{2}) \to \id_L(\frac{d-p-5}{2}) \to \omega_{Y_H}(-1) \to 0$ shows that the divisors of $\omega_{Y_H}(-1)$ are cut on $Y_H$ by surfaces of degree $\delta=\frac{d-p-5}{2}$, containing $L$ but not $T$, residually to $L \cap Y_H$. We may consider surfaces of the form: $H_1 \cup \ldots \cup H_{\delta}$, where $H_1$ contains $L$ and where $H_2,\ldots,H_{\delta}$ are general planes. It follows that our condition is satisfied if $\delta-1 \geq 3$. If $\delta \leq 3$, then $p \geq d-11$ and we conclude as above.\\
If $Q_H$ is the union of two distinct planes, then $Y_H$ is the union of two distinct plane curves. We have: $p_a(Y_H) \geq (\frac{d-p}{2}-1)(\frac{d-p}{2}-2)-1$, because the minimal value for the arithmetical genus of a union of two plane curves of global degree $\delta$ is achieved when each curve has degree $\frac{\delta}{2}$ and if the two components do not intersect. Consequently: $\pi-1 \geq \frac{d^2+p^2-2pd-6d+6p+4}{4}+ \frac{p^2-3p+2}{2}+d-p-r-2$.\\
We may assume that the general hyperplane section of $S$ does not lie on a cubic surface (otherwise $h^0(\id_S(3)) \neq 0$ and $d \leq 3s$), so $\pi-1 \leq \frac{d^2}{8}$. Comparing these two inequalities (and using $r \leq 4$) we obtain: $6p^2-8p-4dp+d^2-4d-32 \leq 0$. If $d \geq 25$ no value of $p$ can satisfy this inequality, so $d \leq 24$ (for all $s$).\epf
\begin{corollary}
\label{Bq=0}
If $dim(\bc_q)=0$, then $r \leq 4$ and $d$ is bounded by a function of $s$. If $s=4$, $d \leq 40$.
\end{corollary}
\textit{Proof:} Since $\bc_q$ is the intersection of the conics $q_H$, $\id_{\bc_q}(2)$ is globally generated, hence $\bc_q$ is contained in a complete intersection of two conics. Recalling that $\rc \subset \bc_q$, it follows that $r \leq 4$ and that $\rc$ does not contain three collinear points. We conclude by lemma \ref{rleq4}.\epf
\begin{lemma}
\label{Bq=1,DnotinS}
Assume $dim(\bc_q)=1$, that $\bc_q$ contains a line $D$ and that $D \not \subset S$. In this case $d \leq s$.
\end{lemma}
\textit{Proof:} Under these assumptions, we claim that the general curve $C$ is smooth. Indeed, let $|L|$ be the linear system cut on $S$ by the hyperplanes containing $D$ and let $B=D \cap S=\{p_1,\ldots,p_r\}$. Clearly $B$ is the base locus of $|L|$ and the general element of $|L|$ is smooth out of $B$. If all curves in $|L|$ were singular at a point $p_i \in B$, it would be $T_{p_i}S \subset H$, $\forall H \supset D$. Anyway the intersection of all $H \supset D$ is nothing but $D$, so this is absurd. The same holds for all $p \in B$. It follows that the singular curves in $|L|$ form a closed subset of $|L|$.\\
Since $D$ is contained in the $\bc_q$, $D$ is $(s-1)$-uple for $\sig$ (see \ref{lem-multiple}). Let $H$ be a general hyperplane through $D$. Then $F=\sig \cap H$ is a degree $s$ surface of $\Pt$ with a line, $D$, of multiplicity $(s-1)$. Such a surface is a projection of a degree $s$ surface $F' \subset \Pq$. We have $S \cap H = C \subset F$ and we may assume $C$ smooth and irreducible. Moreover since $q(S)=0$, C is linearly normal in $\Pt$. Now $C$ is the isomorphic projection of a degree $d$ curve $C' \subset F'$ (in particular $\oc_{C'}(1) \cong \oc_C(1)$). Hence $C'$ is degenerate in $\Pq$ and this implies $d \leq s$.\epf
\begin{lemma}
\label{Bq1=D,DinS}
Assume that the one-dimensional part of $\bc_q$ is a line $D$ and that $D \subset S$. Then $r \leq 1$ and lemma \ref{rleq4} applies.
\end{lemma}
\textit{Proof:} In this case $q_H=D \cup D_H$ and the ${D_H}'s$ are moving. The base locus of the ${D_H}'s$, $\dc$, is either empty or a point, $b$. If $\dc=\emptyset$, then $Y_H \cap \Pi \subset P$ and it follows that $r=0$. Hence we assume from now on that $\dc=\{b\}$.\\
If $b \in D$ we have $\bc_q=D \cup \eta_b$, where $\eta_b$ is the first infinitesimal neighbourhood of $b$ in $\Pi$. Let $x \in Y_H \cap \Pi$ for a general $H$ and let $\xi_x$ be the zero-dimensional subscheme of $Y_H \cap \Pi$ supported at $x$. We will prove the following:\\
\textit{Claim:} Let $x \in Y_H \cap \Pi$, if $\xi_x \not \subset P$ then $x=b$ and, moreover, $\xi_x \subset \eta_b$ if $b \in D$.\\
\textit{Proof of the Claim:} We have $\xi_x \subset S \cap \Pi$. If $\xi_x \not \subset P$ then its residual scheme with respect to $P$ is non empty and so is, a fortiori, the residual scheme of $Z=S \cap \Pi$ with respect to $P$, namely $\rc$. So $\rc$ has a component, $\rc_x$, supported at $x$. Since $\rc \subset \bc_q$, we conclude that $x=b$ or $x \in D$.\\
If $x=b$ and $b \not \in D$, we are over. So we assume $x \in D$. Since $\xi_x \subset q_H$, if $x \neq D \cap D_H$, then $\xi_x \subset D \subset P$: absurd. Thus $x=D \cap D_H$. If $b \in D$ this implies $x=b$ and $\xi_x \subset \eta_b$ (because $\xi_x \subset q_H$). So we may assume $b \not \in D$. In this case the ${D_H}'s$ have no base point on $D$ thus if $H$ is general: $\rc \cap D \cap D_H = \emptyset$: contradiction ($x \in \rc \cap D \cap D_H$).\\
\\
We come back to the proof of the lemma. If $\dc=\{b\}$ and $b \not \in D$ then $Y_H \cap \Pi \subset P$ but for at most one point ($b$), so $Y_HP \geq d-p-1$ and $r \leq 1$.\\
If $\dc=\{b\}$ and $b \in D$, then $\forall x \in Y_H \cap \Pi$, $\xi_x \subset \eta_b$, the residual scheme of $\xi_x$ with respect to $D$ is contained in the residual scheme of $\eta_b$ with respect to $D$, which is $b$. This shows that $Y_HP \geq d-p-1$, hence $r \leq 1$.\epf
\begin{lemma}
\label{Bq=q,qinS}
Assume that $\bc_q$ is a conic $q$ ($q_H=q$ for all $H$). If $q \subset S$, then $r=0$ and lemma
\ref{rleq4} applies.
\end{lemma}
\textit{Proof:} In this case $q \subset P$. Since $Y_H \cap \Pi \subset q_H$, we have $Y_H \cap \Pi \subset P$, hence $Y_HP=d-p$, i.e. $r=0$.\epf
\begin{lemma}
\label{Bq=q,qnotinS}
Assume that $\bc_q$ is a conic $q$ and $q \not \subset S$. Then $d \leq max\{s,20\}$.
\end{lemma}
\textit{Proof:} If no component of $q$ is contained in $S$ (i.e. in $P$), then $Y_H \cap \Pi = Y_H \cap q$ is fixed (otherwise, as $H$ varies, the points of $Y_H \cap \Pi$ will cover a component of $q$). So $Y_H \cap q=\rc$, i.e. $d-p=r$. Since $r=d-2p+P^2$ we get $P^2=p$ and $Y_HP=(H-P)P=0$, this means that $C_H=Y_H \cup P$ is disconnected: absurd.\\
It follows that $q=D \cup L$ with $D \subset S$ and $q \not \subset S$. If $L \neq D$ we have $L \subset \bc_q$, $L \not \subset S$ and we conclude that $d \leq s$ thanks to lemma \ref{Bq=1,DnotinS}.\\
So we may assume $q=2D$, $D \subset P \subset S$ but $2D \not \subset S$ ($2D$ means $D$ doubled in $\Pi$). In this case, for all $H$, $q_H=2D$, so $Q_H$ is tangent to $\Pi$ along $D$. This implies that, for a general $H$, $Q_H$ is either a cone or the union of two distinct planes through $D$. In this latter case $Y_H=P_1 \cup P_2$ and $Y_HD=P_1D+P_2D=d-p$. Since $Y_HD \subset Y_HP$ it follows that $r=0$ and we conclude with lemma \ref{rleq4}.\\
From now on we assume that for a general $H$, $Q_H$ is a cone and $D$ a ruling of $Q_H$. If $d-p$ is even, $Y_H$ is a complete intersection $(\frac{d-p}{2},2)$, then $p_a(Y_H)=\frac{d^2-2pd-4d+p^2+4p+4}{4}$ and so $\pi-1=\frac{d^2-2pd-4d+p^2+4p+4}{4} + \frac{p^2-3p+2}{2}+d-p-r-2$. Now $Y_H \cap D \subset Y_H \cap P$, then $Y_HD=\frac{d-p}{2} \leq d-p-r=Y_HP$, i.e. $r \leq \frac{d-p}{2}$ and it follows that $\pi-1 \geq \frac{d^2-2pd-4d+p^2+4p+4}{4} + \frac{p^2-3p+2}{2}+ \frac{d-p}{2}-2$. Now comparing this expression with $\pi-1 \leq \frac{d^2}{8}$ (we can suppose as usual $h^0(\id_C(3))=0$) we get: $6p^2-8p-4dp+d^2-4d \leq 0$. If $d \geq 21$ there are no values of $p$ satisfying the inequality, then $d \leq 20$.\\
If $d-p$ is odd, $Y_H$ is linked to a line by a complete intersection $(\frac{d-p+1}{2},2)$ and it turns out $p_a(Y_H)=\frac{d^2-2dp+p^2-4d+4p+3}{4}$. Since $Y_HD=\frac{d-p+1}{2} \leq Y_HP=d-p-r$ we have $r \leq \frac{d-p-1}{2}$. Hence we can write $\pi-1 \geq \frac{d^2-2dp+p^2-4d+4p+3}{4}+ \frac{p^2-3p+2}{2}+\frac{d-p+1}{2}-2$. If we compare this with $\pi-1 \leq \frac{d^2}{8}$ and arguing as before we obtain $d \leq 20$.\epf\\
\\
The proof of \ref{th1} and \ref{th2} follows from \ref{small-genus}, \ref{rleq4}, \ref{Bq=0}, \ref{Bq=1,DnotinS}, \ref{Bq1=D,DinS}, \ref{Bq=q,qinS}, \ref{Bq=q,qnotinS}.
\begin{remark}
\label{general}
Actually we believe that there are very few smooth surfaces on such hypersurfaces. For example consider the following situation:\\
Assume that the blowing-up of $\Pi$, $\tilde{\Sigma}\to \Sigma$, yields a desingularization of $\Sigma$, so we have a double covering $T \to \Pi$ and $\tilde{S}$ mapping to $S$. Since $T$ and $\tilde{S}$ are two divisors on the smooth threefold $\tilde{\Sigma}$, if they intersect, they intersect along a curve. We conclude that $S \cap \Pi = P$ and all the points of $Y_H \cap \Pi$ lie on $P$.\\
Now assume that for general $H$, $Q_H$ is a smooth quadric. Observe that the $Q_H$ are parametrized by a smooth rational curve ($\simeq \Pun$). Let $\pc$ denote the curve parametrizing the rulings of the quadrics $Q_H$. We get a degree two covering $f:\pc \to \Pun$ which is ramified at the points corresponding to singular $Q_H$. Assume $\pc$ is irreducible. With this assumption the curve $Y_H \subset Q_H$ has bidegree $(a,a)$ (otherwise following the $a$ ruling would yield a section of the covering, which is impossible since $g(\pc )>0$ because $f$ is ramified in more than two points).\\
Now consider the exact sequence of residuation with respect to $\Pi$:\\
$$0 \to \id _{Y_H}(-1) \to \id _C \to \id _{P,\Pi} \to 0$$
Since $Y_H$ is a.C.M., it follows that $C = Y_H \cup P$ is a.C.M. too. Hence $S$ is a.C.M. and $h^0(\id _S(3)) \geq h^0(\id _C(3))\neq 0$. This implies $d(S) \leq 3s$. (Notice that we didn't assume $q(S)=0$.) Observe that the assumption that $S$ is smooth is necessary in order to apply Lemma \ref{lem1} and to conclude that $C = Y_H\cup P$ with $Y_H \subset Q_H$.
\end{remark}
\begin{remark}
There exist integral hypersurfaces in $\Pq$ such that the degree of the smooth surfaces contained in them is bounded. Indeed it is enough to take a non linearly normal hypersurface in $\Pq$, recalling that the only non linearly normal smooth surface in $\Pq$ is the Veronese. The simplest example is the Segre cubic hypersurface. The previous results seem to indicate that this behaviour can happen also on some linearly normal hypersurfaces. From a "codimension two" point of view this is in contrast with the following proposition.
\end{remark}
\begin{proposition}
Let $S \subset \Pt$ be an integral surface, then $S$ contains smooth curves of arbitrarily high degree.
\end{proposition}
\textit{Proof:} If $S$ has singular locus of dimension $\leq 0$, this follows from Bertini. If $Sing(S)$ has dimension 1, we consider the normalization $p:\tilde{S} \to S$ of $S$, then $dim(Sing(\tilde{S}))\leq 0$. Let $C$ be the non-normal locus in $S$, $D=p^{-1}(C)$. Let $\delta$ be a very ample linear system on $\tilde{S}$. The general $X \in \delta$ is smooth and doesn't pass through any singular point of $\tilde{S}$. We want to show that for $X \in \delta$ general, $p_|:X \to S$ is an embedding.
Since $p$ is an isomorphism outside $D$, we only have to consider the points in $X \cap D$. Let $x \in C$, the curves of $\delta$ passing through two points of $p^{-1}(x)$ form a subspace of codimension $2$. Letting $x$ vary in $C$, we see that the curves of $\delta$ intersecting a fibre $p^{-1}(x)$ in more than one point constitute a subspace of codimension $\geq 1$, hence for general $X \in \delta$, $p_|:X \to S$ is injective.\\
Since there are only finitely many points where $dp$ has rank zero, we may assume that for $y\in D$, $dp_y:T_y\tilde{S} \to T_{p(y)}S$ has rank one. The curves of $\delta$ passing through $y$ and having tangent direction $Ker(dp_y)$ at $y$ form a subspace of codimension $2$ of $\delta$. Letting $y$ vary in $D$ we get a subspace of codimension $1$. So for general $X \in \delta$, $dp_|$ is everywhere injective.\epf

\bigskip

\noindent
Address of the authors:\\
Dipartimento di Matematica\\
via Machiavelli, 35\\
44100 Ferrara (Italy)\\
Email: phe@dns.unife.it (Ph.E.), chiara@dm.unife.it (C.F.)

\end{document}